\documentclass{sn-jnl}% Math and Physical Sciences Numbered Reference Style 
\usepackage[numbers]{natbib}
\usepackage{graphicx}%
\usepackage{multirow}%
\usepackage{amsmath,amssymb,amsfonts}%
\usepackage{amsthm}%
\usepackage{mathrsfs}%
\usepackage{dsfont}%
\usepackage[title]{appendix}%
\usepackage{xcolor}%
\usepackage{textcomp}%
\usepackage{manyfoot}%
\usepackage{booktabs}%
\usepackage{algorithm}%
\usepackage{algorithmicx}%
\usepackage{algpseudocode}%
\usepackage{listings}%
\usepackage{todonotes}%
\usepackage{mathtools} % for \coloneqq
\usepackage{hyperref} % For clickable links in the caption
\usepackage{xspace}
\usepackage{caption}
\usepackage{subcaption}
\usepackage{tablefootnote}
\usepackage{orcidlink}

\newcommand{\win}[1]{\textbf{\textcolor{blue}{#1}}}
\newcommand{\loss}[1]{\textit{\textcolor{red}{#1}}}

\setuptodonotes{inline, color=blue!30}

\theoremstyle{thmstyleone}%
\newcommand{\solver}[1]{\textsc{#1}}

\newcommand{\scip}{\solver{{\sc Scip}}\xspace}
\newcommand{\papilo}{\solver{Papilo}\xspace}
\newcommand{\scipjack}{\solver{{\sc Scip-Jack}}\xspace}
\newcommand{\scipsdp}{\solver{{\sc Scip-sdp}}\xspace}

\newcommand{\fscip}{\solver{FiberScip}\xspace}

\newcommand{\abscon}{\solver{absconpseudo}\xspace}
\newcommand{\exact}{\solver{exact}\xspace}
\newcommand{\galena}{\solver{galena}\xspace}
\newcommand{\minisat}{\solver{minisat+}\xspace}
\newcommand{\naps}{\solver{naps}\xspace}
\newcommand{\openwbo}{\solver{openwbo}\xspace}
\newcommand{\parlspbo}{\solver{parls-pbo}\xspace}
\newcommand{\prs}{\solver{prs}\xspace}
\newcommand{\pueblo}{\solver{pueblo}\xspace}
\newcommand{\roundingsat}{\solver{roundingsat}\xspace}

\newcommand{\SatforjPBRes}{\solver{sat4j pb res}\xspace}
\newcommand{\SatforjPBCP}{\solver{sat4j pb res}\xspace}

\newcommand{\Z}{\mathds{Z}}
\newcommand{\N}{\mathds{N}}

\usepackage{todonotes,xspace}

\begin{document}

\title[Pseudo-Boolean Solving with SCIP]{State-of-the-art Methods for Pseudo-Boolean Solving with SCIP}

\author*[1]{\fnm{Gioni} \sur{Mexi} \orcidlink{0000-0003-0964-9802}}\email{mexi@zib.de}% OrcID 0000-0003-0964-9802
\author[3]{\fnm{Dominik} \sur{Kamp}\orcidlink{0009-0005-5577-9992}\email{dominik.kamp@uni-bayreuth.de}}% OrcID 0009-0005-5577-9992

\author[1]{\fnm{Yuji} \sur{Shinano}\orcidlink{0000-0002-2902-882X}}\email{shinano@zib.de}% OrcID 0000-0002-2902-882X
\author[1]{\fnm{Shanwen} \sur{Pu}\orcidlink{0009-0005-5034-8341}}\email{pu@zib.de}% OrcID 0009-0005-5034-8341
\author[1]{\fnm{Alexander} \sur{Hoen}\orcidlink{0000-0003-1065-1651}}\email{hoen@zib.de}% OrcID 0000-0003-1065-1651
\author[1]{\fnm{Ksenia} \sur{Bestuzheva}\orcidlink{0000-0002-7018-7099}}\email{bestuzheva@zib.de}% OrcID 0000-0002-7018-7099
\author[4]{\fnm{Christopher} \sur{Hojny}\orcidlink{0000-0002-5324-8996}}\email{c.hojny@tue.nl}% OrcID 0000-0002-5324-8996
\author[5]{\fnm{Matthias} \sur{Walter}\orcidlink{0000-0002-6615-5983}}\email{matthias@matthiaswalter.org}% OrcID 0000-0002-6615-5983
\author[6]{\fnm{Marc E.} \sur{Pfetsch}\orcidlink{0000-0002-0947-7193}}\email{pfetsch@mathematik.tu-darmstadt.de} % OrcID 0000-0002-0947-7193
\author[1,2]{\fnm{Sebastian} \sur{Pokutta}\orcidlink{0000-0001-7365-3000}}\email{pokutta@zib.de} % OrcID 0000-0001-7365-3000
\author[1,2]{\fnm{Thorsten} \sur{Koch}\orcidlink{0000-0002-1967-0077}}\email{koch@zib.de} % OrcID 0000-0002-1967-0077

%----------------------------------------------------------------
\affil[1]{\orgname{Zuse Institute Berlin}, \orgaddress{\country{Germany}}}

\affil[2]{\orgname{TU Berlin}, \orgaddress{\country{Germany}}}

\affil[3]{\orgname{University of Bayreuth}, \orgaddress{ \country{Germany}}}

\affil[4]{\orgname{Eindhoven University of Technology}, \orgaddress{ \country{Netherlands}}}

\affil[5]{\orgname{University of Twente}, \orgaddress{ \country{Netherlands}}}

\affil[6]{\orgname{TU Darmstadt, Department of Mathematics}, \orgaddress{\country{Germany}}}

% \affil[2]{\orgdiv{Department}, \orgname{Organization}, \orgaddress{\street{Street}, \city{City}, \postcode{10587}, \state{State}, \country{Country}}}

\abstract{
The Pseudo-Boolean problem deals with linear or polynomial constraints with integer coefficients over Boolean variables. The objective lies in optimizing a linear objective function, or finding a feasible solution, or finding a solution that satisfies as many constraints as possible.
In the 2024 Pseudo-Boolean competition, solvers incorporating the \scip framework won five out of six categories it was competing in. %\scip and \fscip won four out of six categories. 
%Furthermore, a solver incorporating \scip won another category.
%
From a total of 1,207 instances, \scip successfully solved 759, while its parallel version \fscip solved 776.
Based on the results from the competition, we further enhanced \scip's Pseudo-Boolean capabilities.
This article discusses the results and presents the winning algorithmic ideas.

}




\maketitle

\section{Introduction}
\label{sec:intro}

%For a quick overview look at YT: https://www.youtube.com/watch?app=desktop&v=yvXy7Nzn-IA

The \emph{Boolean Satisfiability problem} (SAT) is one of the most studied and fascinating problems of computer science~\cite{FrancoMartin2009HistoryofSAT}.
Modern SAT solvers are typically based on the \emph{conflict driven clause learning} (CDCL) paradigm introduced by \cite{MarquesSakallah1999Grasp}.
Over recent years, these solvers have significantly improved their performance, allowing them to solve even large-scale problems.
Therefore, people searched for new fields of applications and applied SAT solvers to \emph{Pseudo-Boolean (PB) problems}, i.e., integer (linear) constraints over Boolean variables. 

One way to solve these problems is the transformation to \emph{conjunctive normal form} (CNF), allowing CDCL-based solvers to solve the problem. 
This approach is, for example, used by solvers such as, e.g., \minisat\cite{En2006TranslatingPC}, \openwbo\cite{MartinsManquinho2014OpenWbo,JoshiMartins2015OpenWbo}, and \naps\cite{SakaiNabeshima2015Construction}.
A similar strategy is employed by one of the algorithms in \SatforjPBRes\cite{LeBerreParrain2010Sat4J22} that keeps the original data input but derives new information only in the clause format.

A different approach is to adapt/generalize the methods from SAT to Pseudo-Boolean: 
The solvers \prs\cite{DixonGinsberg02InfernceMethods}, \galena\cite{ChaiKuehlmann05Galena}, \pueblo~\cite{SheiniSakallah06Pueblo} and \SatforjPBCP~\cite{LeBerreParrain2010Sat4J22} generalize the cutting-plane proof system~\cite{Hooker92GeneralizedResolution} based on saturation to Pseudo-Boolean, while  \roundingsat~\cite{elfers18DivideAndConquer,DevriendtGochtDemirovicNordstrom2021CuttingtotheCoreofPBO} and its fork \exact{} employ a conflict-driven search, but use division or rounding instead of saturation to keep coefficients small.
Further approaches are applying constraint programming techniques used in \abscon~\cite{HemeryLecoutre2006AbsconPseudo} or even local search strategies, for example, \parlspbo~\cite{ChenLinHuCai2024ParLSPBO}.
For a more complete overview of PB-solvers, see \cite[Chapter 2]{Biere2009HandbookofSatisfiability}.

Pseudo-Boolean problems can also be formulated as \emph{integer linear programs} (ILP).
This allows (mixed) integer-programming solvers to solve PB problems.
In the past PB competitions, \scip\cite{Achterberg07Thesis,SCIP9} (see also Section \ref{sec:scip}) consistently demonstrated top performance in nonlinear and optimization categories.
Throughout all four competitions~\cite{PB09,PB10,PB11,PB12} it previously participated in, \scip won first or second places on the benchmarks with nonlinear decision problems and linear optimization problems and the first places on the benchmark with nonlinear optimization problems.
It gradually lost its position in the ranking for linear decision problems, moving from second place in 2009 to nineteenth place in 2012, but strengthened its position in the rankings for Weighted Boolean Optimization (WBO), moving from fifth place in 2009 to first place in 2012 both in categories of problems with all soft clauses and problems with mixed clauses.

% Results:
% https://www.cril.univ-artois.fr/PB09/results/ranking.php?idev=28
% https://www.cril.univ-artois.fr/PB10/results/ranking.php?idev=36
% https://www.cril.univ-artois.fr/PB11/results/ranking.php?idev=54
% https://www.cril.univ-artois.fr/PB12/results/ranking.php?idev=67

This paper aims to showcase the viability of \scip and its multithreaded version \fscip for solving various types of PB problems. After defining the problem and providing a brief overview of the \scip and \fscip frameworks in Section \ref{sec:intro}, we highlight key modifications and new features implemented in \scip to enhance its handling of Pseudo-Boolean problems in Section \ref{sec:features}. 
% In Section \ref{sec:fscip}, we demonstrate the simplicity of using \fscip in racing mode, showcasing how minimal customization can enable parallelization for Pseudo-Boolean problems. 
In Section \ref{sec:results}, we present the results of the PB24 competition across various problem categories, followed by a discussion of post-competition improvements.

\subsection{Problem Definition}
\label{sec:definition}

Pseudo-Boolean problems involve $n$ Boolean variables\footnote{The normalized form used in PB literature to express constraints typically contains literals. 
    A \emph{literal} of a variable $x$ is either the variable itself or its negation $\neg x= 1-x$.
    Here, we use the more general definition of IP.
}.
A \emph{linear Pseudo-Boolean constraint} over $n$ binary variables $x_1, \dots, x_n$ is defined as
\begin{align*}
    \label{eq::nonlinearPbconstraint}
    \sum_{j=1}^n a_j x_j \geq b
\end{align*}
with $a\in \Z^n$ and $b\in \Z$.
A \emph{non-linear Pseudo-Boolean constraint} over $n$ binary variables is defined as 
\begin{align}
    \sum_{i=1}^t a_i \prod_{j\in M_i} x_{j}\geq b
\end{align}
with $t\in \N$, $a\in \Z^t$, $b\in \Z$ and $M_i \subseteq \{1,\dots,n\}$ for all $i\in \{1,\dots,t\}$.
Each monomial $\prod_{j\in M_i} x_{j}$ can be reformulated as an AND constraint by introducing new binary variables
\begin{align*}
    z_i = \bigwedge_{j\in M_i} x_{j}.
\end{align*}
The new variables $z_i$ can be inserted into \eqref{eq::nonlinearPbconstraint}, resulting in a linear constraint $a^T z \geq b$. 
It is well-known that AND constraints can be linearized by linear equations as discussed in \cite{berthold2008solving}.

Hence, PB problems can be reformulated as pure \emph{integer linear program} (ILP) over  $n$ binary variables in the following form~\cite{berthold2008solving}:
\begin{align}
    \min \;&c^T x \nonumber\\
     &Ax \geq b, \label{eq:PB}\\
     &x \in \{0, 1\}^n \nonumber
\end{align}
with $A\in \Z^{m\times n}$, $b\in \Z^{m}$ and  $c\in \Z^{n}$.
If the problem does not have an objective, it is called a \emph{decision} problem; otherwise, an \emph{optimization} problem.

The Pseudo-Boolean competition also contains problems in which some or all constraints can be relaxed, and the number of satisfied constraints is maximized (or the number of unsatisfied constraints minimized).
Such constraints are referred to as \emph{soft}. In contrast, constraints that cannot be relaxed are called \emph{hard}.
Soft constraints can be modeled using additional binary variables $y$ for each such constraint. Then, a so-called indicator constraint $y = 1 \to \sum_{j = 1}^n a_j x_j \geq b$ is used, meaning the constraint $\sum_{j = 1}^n a_j x_j \geq b$ has to hold if $y = 1$. In the case in which all constraints are allowed to be relaxed, \eqref{eq:PB} turns into
\begin{align}\label{eq:PBindicator}
    \max \;& \sum_{i=1}^m y_i \nonumber\\
     & y_i = 1 \to \sum_{j = 1}^n a_{ij} x_j \geq b_i \quad \forall\, i = 1, \dots, m, \\
     &x \in \{0, 1\}^n,\; y \in \{0,1\}^m, \nonumber
\end{align}
where $A = (a_{ij}) \in \Z^{m \times n}$.

\subsection{SCIP and FiberSCIP}
\label{sec:scip}
\scip \cite{Achterberg07Thesis, SCIP9} is an open-source Constraint Integer Programming (CIP) solver.
The development of the concept of CIP was motivated by the fact that the majority of solvers for SAT, MIP, and CP work in the spirit of branch-and-bound, which means that they recursively subdivide the problem instance, yielding a so-called branch-and-bound-tree, whose nodes represent subproblems of the original instance.
CIPs generalize these problem classes and are defined as finite-dimensional optimization problems with arbitrary constraints and a linear objective function that satisfies the following property: If all integer variables are fixed, the remaining subproblem must form a linear or nonlinear program (LP or NLP).

As a CIP solver, \scip seamlessly combines techniques from SAT, MIP, and CP and offers greater flexibility in how problems are defined and the potential for extensions to further problem classes.
PB problems are one natural application of such a framework since they are a special case of CIP.

\scip employs an LP-based branch-and-cut algorithm, enhanced with numerous additional techniques crucial for the solver's performance. In the context of PB problems, the following are especially important:
\begin{itemize}
    \item \textit{Constraint propagation} analyzes the constraints of the current subproblem and the local domains of the variables to infer additional valid constraints and domain reductions, thereby restricting the search space.
    \item \textit{Conflict analysis} enables the solver to learn from infeasible subproblems by deducing short, globally valid conflict clauses from a series of branching decisions that led to infeasibility. These clauses enable the solver to prune other parts of the search tree and apply non-chronological backtracking.
    \item \textit{Cutting plane techniques} tighten the LP-relaxation of the problem by adding additional linear inequalities that cut off a possibly fractional solution of the LP relaxation.
    \item \textit{Primal heuristics} aim to produce feasible solutions of high quality.
    \item \textit{Presolving procedures} aim at transforming the original problem into an equivalent problem that is easier to solve.
    \item \textit{Restarts} abort the search process if a certain amount of global problem reductions has been triggered in the early steps and restart the search from scratch. The motivation is to use the knowledge obtained in previous runs by reapplying other presolving mechanisms to the reduced problem, as well as procedures that are only applied at the root node.
    \item \textit{Indicator Constraints}
    are handled by branching on the decision variables $y_i$ in~\eqref{eq:PBindicator}. If \scip determines that the slack of the corresponding constraints is small enough, it adds a linear inequality to the LP using big-M values, i.e., for the indicator constraint $y = 1 \to \sum_{j=1}^n a_j x_j \geq b$ it would add $\sum_{j=1}^m a_j x_j + M\, (1 - y) \geq b$, where $M > 0$ is large enough such that $\sum_{j=1}^m a_j x_j + M \geq b$ holds for all feasible solutions $x$.
    Restarting is one effective feature for Pseudo-Boolean Problems, especially if all $y_i$ variables are fixed. After the restart, one obtains a ``regular'' PB problem, for which more presolving can be performed. 
\end{itemize}

In 2024, the parallel solver \fscip~\cite{shinano2018fiberscip} participated for the first time in the PB competition alongside sequential \scip.
\fscip is a shared memory parallelization of \scip created using the UG framework~\cite{Shinano2019}, which provides a systematic way to parallelize branch-and-bound solvers.

\fscip is built by parallelizing the tree search from ``outside'' of \scip by dynamically splitting up the search tree and maintaining each subtree in a separate \scip instance.
\scip's tree search is performed in each subtree, benefiting directly from the state-of-the-art solving techniques of \scip. Only the distribution of subproblems to \scip solvers is controlled from the ``outside'' to balance the workload effectively.

\section{Added Features for Pseudo-Boolean Problems\textbf{}}
\label{sec:features}
In this section, we highlight key modifications and new features implemented in \scip to handle Pseudo-Boolean problems better, enhancing both the solver's robustness and efficiency for this domain.

\subsection{RLT Cuts for AND Constraints}

The reformulation-linearization technique (RLT) was first proposed by Adams and
Sherali~\cite{adams1986tight,adams1990linearization,adams1993mixed} for linearly constrained
bilinear problems with binary variables and can be applied to mixed-integer, quadratic, and general polynomial problems. RLT constructs valid polynomial constraints and then linearizes these constraints by using relations given in the problem when possible and applying relaxations otherwise.

\scip{} generates RLT cutting planes for bilinear terms, including bilinear terms participating in nonlinear relations implicitly given by mixed-integer linear constraints.
Consider a linear constraint:
$\sum_{j=1}^n a_{ij}x_j \geq b_i.$
Multiplying this constraint by nonnegative bound factors $(x_k-\underline{x}_k)$ and $(\overline{x}_k-x_k)$, where, in the case of binary variables considered in this paper, $\underline{x}_k = 0$ and $\overline{x}_k = 1$, yields valid nonlinear inequalities.
For example, multiplication by the lower bound factor (that is, the variable $x_k$ itself) yields:
\[
\sum_{i=1}^n a_{ij}\, x_j\,x_k \geq b_i\,x_k.
\]
The above nonlinear inequality is then linearized to obtain a valid linear inequality.
This linearization step substitutes products with linear expressions if a relation between a product term and a linear expression is given or detected in the model, with equivalent expressions if cliques enable such reformulations, and applies McCormick relaxations otherwise.
Square terms $x_j^2$ are substitutes with $x_j$ since $x_j^2 = x_j$ holds for a binary variable.

An efficient separation algorithm prevents excessive time spent in RLT cut separation and is crucial for efficiently applying RLT in practice. The implicit product detection and separation algorithm was introduced by ~\cite{AchterbergBestuzhevaGleixner2024}.

Since in Pseudo-Boolean instances, the products of variables are represented by AND constraints, the RLT cut separator had to be extended to search for products in these constraints in addition to nonlinear and linear constraints.

\subsection{Flower Inequalities}
\label{sec:mods:flower}

A new separator dealing with so-called \emph{multilinear} constraints was added to \scip. These are constraints of the form $z = x_1 \cdot x_2 \cdot \dotsc \cdot x_k$ for binary variables $x_i \in \{0,1\}$.
Each such constraint was already treated individually in \scip as it models the logical AND constraint.
The purpose of the new separator is to generate cutting planes that strengthen the LP relaxation in the presence of multiple such constraints.
The collection of all constraints can be represented via a hypergraph $G = (V,E)$, which has a node $v \in V$ per variable $x_v$ that appears in at least one AND constraint and a hyperedge $e \in E$ per constraint, where the product variable satisfies $z_e = \prod_{v \in e} x_e$.

This hypergraph is constructed after presolving the problem by inspecting all AND constraints and scanning expression trees of nonlinear constraints.
Moreover, \emph{overlap sets}, which are sets of the form $e \cap f$ for $e$, $f \in E$, are gathered using hash tables, from which \emph{compressed sparse row} and \emph{-column} representations of the incidences among hyperedges, overlap sets and nodes are computed.

While there exist many classes of valid inequalities~\cite{CramaR17,DelPiaK17,DelPiaK18,DelPiaD21,DelPiaW22}, we restricted the current implementation to \emph{$k$-flower inequalities}~\cite{DelPiaK18} for $k=1,2$, since the separation problem can be solved very efficiently once the overlap sets are available.
Such an inequality is parameterized by one \emph{center} edge $e \in E$ and $k$ \emph{neighbor edges} $f_1,f_2,\dotsc,f_k \in E$, all of which must intersect $e$.
It reads
\[
  z_e + \sum_{i=1}^k (1-z_{f_i}) + \sum_{v \in R} (1-z_v) \geq 1,
\]
where $R \coloneqq e \setminus \bigcup_{i=1}^k f_i$ denotes the nodes of $e$ that are covered by no neighbor edge.
By simple enumeration, the separation problem can easily be solved in time $\mathcal{O}(|E|^{k+1})$.
Under the reasonable assumption that the size $|e|$ of every edge is bounded by a constant, this can be reduced to $\mathcal{O}(|E|^2)$ as done in the implementation~\cite{DelPiaKS20}.
By exploiting overlap sets using the implemented data structures, this can be reduced to $\mathcal{O}(|E|)$.

\subsection{Symmetries}

A symmetry of a Pseudo-Boolean problem with variables~$x \in \{0,1\}^n$ is a map $\varphi\colon \{0,1\}^n \to \{0,1\}^n$ satisfying, for every~$\bar{x} \in \{0,1\}^n$, the following two properties:
\begin{itemize}
\item $\varphi(\bar{x})$ is a feasible solution if and only if~$\bar{x}$ is feasible, i.e., $\varphi$ preserves feasibility;
\item both~$\varphi(\bar{x})$ and~$\bar{x}$ have the same objective value, i.e., $\varphi$ preserves the objective.
\end{itemize}
It is well-known that handling symmetries can greatly improve the performance of SAT, CP, and MIP solvers because the unnecessary exploration of symmetric parts of the search space can be avoided.
Since detecting all symmetries of a Pseudo-Boolean problem is NP-hard~\cite{Margot2010}, \scip only detects symmetries that keep the problem syntax invariant~\cite{SCIP5}.

The symmetries that can be handled by \scip are permutation symmetries, which correspond to reorderings of the variable vector~$x$, reflection symmetries, which map some variable~$x_i$ to its binary negation~$1 - x_i$, and combinations thereof.
To detect such symmetries, \scip generates so-called symmetry detection graphs~\cite{Hojny2024} for each constraint, which is a colored graph, all of whose automorphisms correspond to symmetries of the constraint.
In particular, symmetry detection graphs have been implemented for linear and Pseudo-Boolean constraints.
By combining the symmetry detection graphs of all constraints in a suitable way, symmetries of an entire Pseudo-Boolean problem can be detected; cf.~\cite{Hojny2024}.

\scip offers a variety of different symmetry handling methods.
For Pseudo-Boolean problems, we observed that handling both permutation and reflection symmetries is beneficial.
This is in contrast to general mixed-integer linear problems, where reflection symmetries seem to occur less often~\cite{Hojny2024}.
After computing permutation and reflection symmetries, \scip tries to detect some structure of the corresponding symmetry group heuristically:
It checks whether some of the problem's variables can be arranged in a matrix~$X \in \{0,1\}^{p \times q}$ such that the symmetry group acts on the matrix~$X$ by permuting its columns; we refer to such a symmetry as an \emph{orbitopal symmetry}; cf.~\cite{KaibelPfetsch2008}.
Moreover, it checks whether a reflection symmetry exists that simultaneously reflects the variables of an entire column of~$X$; we refer to such a symmetry as \emph{reflective orbitopal symmetry}.

To handle reflective orbitopal symmetries for matrix~$X \in \{0,1\}^{p \times q}$, \scip adds the inequalities~$X_{1,1} \geq X_{1,2} \geq \dots \geq X_{1,q} \geq 0$ and enforces that the columns of~$X$ are sorted lexicographically.
The latter is achieved by the  \emph{orbitopal fixing}~\cite{BendottiEtAl2021} propagation algorithm.
For orbitopal symmetries without reflections, \scip does not add the linear inequalities used for reflective orbitopal symmetries.
Instead, \scip checks whether the constraints of a Pseudo-Boolean problem guarantee that, for every solution, each row of~$X$ has at most or exactly one~1-entry.
In this case, orbitopal symmetries for~$X$ are also handled by enforcing the columns of~$X$ to be sorted lexicographically via a specialized version of orbitopal fixing that takes bounds on the number of~1-entries into account~\cite{KaibelEtAl2011} (chronologically, the orbitopal fixing algorithm by~\cite{KaibelEtAl2011} has been developed first though).
If no such bound on the number of 1-entries is detected, a variant of orbitopal fixing is used in which the position of the problem's variable inside the matrix~$X$ is not fixed a priori.
Instead, the position is dynamically decided during the solving process based on the branching history~\cite{DoornmalenHojny2024a}.

If \scip does not detect orbitopal symmetries, the following two mechanisms are applied simultaneously to handle symmetries.
The first mechanism is the propagation algorithm orbital fixing~\cite{OstrowskiEtAl2011}, which detects variables fixed to~0 or~1 
%by branching decisions or other reasons 
and fixes some symmetric variables to the same value.
The second mechanism takes a single permutation or reflection symmetry~$\varphi$ and derives variable fixings from a symmetry handling constraint that enforces a solution~$x$ to be lexicographically not smaller than its image~$\varphi(x)$.
Note that, to be compatible with orbital fixing, this lexicographic order needs to be adapted to the branching history; see~\cite{DoornmalenHojny2024a} for details.

\subsection{Numerical Solution Handling}

During preparation for the competition, we had to overcome some numerical challenges crucial to the correctness of the resulting solutions. For example, on the linear decision subset~\texttt{submitted-PB06/manquiho/Aardal\_1/normalized-prob3} of the PB competition 2006 derived from the set~\cite{Aardal2002} of hard equality constraints, the default setup provided a solution although the instance is known to be infeasible. This feasibility problem consists of a single equality constraint with $52$ variables and the right-hand side $b = 5842800$. The suggested solution is strictly binary but has a constraint activity of $b - 2$
which is within the acceptable range for a relative feasibility tolerance of $\epsilon_f\coloneqq 10^{-6}$. Therefore, we decided to instead apply an absolute tolerance for all relevant feasibility checks, thereby ensuring that our solutions fully complied with the competition's criteria for exactness, regardless of the magnitude of the involved constraints.

Another numerical artifact arises from using floating-point values for the binary variables, which are typically rounded to exact binaries only after the solving process is finished. The linear equality in instance~\texttt{normalized-prob3} has a coefficient $5567264$. Therefore, even when applying the absolute tolerance of $\epsilon_f$, the corresponding value of the exact integral solution above can be increased by a sufficient approximation of $2/5567264 = 1/2783632 \leq \epsilon_f$ to satisfy the constraint. In general, let $C$ denote the coefficient row vector of some constraint. Given an arbitrary tolerable solution, the maximal activity change, when rounding this solution, is given by $\lVert C \rVert_, \epsilon_f$. Since Pseudo-Boolean constraints are formulated with exact integers by design, the activity of the considered tolerable solution is within $\epsilon_f$ around the exact integral side, and the corresponding rounded solution reaches an exact integral activity so that the constraint can only become violated if the absolute activity change is at least $1 - \epsilon_f$. Hence, the rounded solution satisfies this constraint if
\[
\epsilon_f < \frac{1}{\lVert C \rVert_1 + 1}
\]
holds.
The \emph{intsize}, as defined by the OPB format specification, is the number of bits required to represent, for any constraint, the sum of the absolute value of all integers that appear in the constraint.
For an intsize $s$, this exclusive upper bound is at least $2^{s}$. Consequently, for every absolute tolerance $\epsilon_f < 2^{-s}$, a tolerable solution will be reliably rounded to an exactly feasible solution, provided the activity evaluation is sufficiently precise. This issue also occurs on several instances of the new linear decision subset~\texttt{devriendt/bitvector\_multiplication\_selection/equalities/array\_comm} even for small intsizes, which this approach could handle. However, adapting the general tolerance to below $2^{-s}$ can become very restrictive already for moderate intsizes $s$ up to the internal limit $49$ set in \scip. For this reason, to avoid unnecessary rejections of exactly feasible solutions, we decided to keep the tolerance value and instead round every solution before checking feasibility. In this way, the solving process proceeds if the rounded counterpart of a solution turns out to be infeasible.
As a result, using an absolute tolerance and applying the pre-check polishing made it possible for \fscip with 20 cores to solve the hard instance~\texttt{normalized-prob3} correctly to proven infeasibility.

\subsection{Large Integer Heuristics}
% \todo{Shanwen}

To accommodate problems with intsize larger than 49, we modified the Feasibility Jump  heuristic~\cite{luteberget2023feasibility} to support purely integer operations, eliminating reliance on floating-point arithmetic. This adjustment is crucial for maintaining precision and efficiency, particularly when handling large integer values, where floating-point calculations can introduce overflows and computational inefficiencies. Our implementation supports GMP (GNU Multiple Precision)~\cite{gmplib} integers for scenarios requiring handling very large integer values. However, the standard 64-bit `long' type suffices for the current Pseudo-Boolean competition.

In Pseudo-Boolean problems, the computation of the jump value is notably simplified, as it only involves comparing the constraint violations for $x_j=0$ and $x_j=1$. This binary nature allows for rapid evaluation of which assignment minimizes constraint violations and enhances the heuristic’s performance. Furthermore, in this integer-only version, the weights and the scaling factors for constraints are redefined using integer values, ensuring all arithmetic operations remain within the integer domain.
In the benchmark set for the competition, there were, in total, only four instances with an intsize exceeding 49, and our special Feasibility Jump heuristic was able to find a solution for one of them.

\subsection{Pseudo-Boolean Cut-based Conflict Analysis}
\label{sec:conflicts}

Following the competition, we extended \scip to include the conflict analysis technique introduced by~\cite{mexi2024cutbasedconflictanalysismixed}, which generalizes Pseudo-Boolean conflict analysis to MILP. This approach interprets conflict analysis as a sequence of linear combinations, integer roundings, and cut generation. In contrast to traditional SAT-based conflict analysis used in MIP solvers, which operates on clauses extracted from constraints rather than the constraints themselves, this method leverages the full expressive power of the linear constraints. For further details, we refer the reader to \cite{MexiBertholdGleixneretal.2023, mexi2024cutbasedconflictanalysismixed}.

Our analysis in Section~\ref{sec:postcompetition} demonstrates its potential to enhance solver performance.
Specifically, by re-evaluating the competition settings with this method enabled, \scip is able to solve eight additional instances within the time limit and improve the average running time by 2\% on instances where such conflict constraints are generated.

\section{Competition Results}
\label{sec:results}

This section presents a re-evaluation of the results of the PB24 competition, highlighting the performance of \scip and \fscip across various problem categories.

\subsection{Categories}
% \AH{maybe we should pick up the definitions of the introduction. NLC: at least on constraints of type \ref{eq::nonlinearPbconstraint}, DEC: a PBS instance}
The competition problems are categorized based on two main criteria: the linearity of constraints and the existence of an objective function. Additionally, problems are further classified based on the presence of hard or soft clauses. The categories used in this competition are summarized as:
\begin{itemize}
    \item \textit{\textsc{opt-lin}}: Optimization problems with all linear constraints.
    \item \textit{\textsc{dec-lin}}: Decision problems with all linear constraints.
    \item \textit{\textsc{opt-nlc}}: Optimization problems with at least one nonlinear constraint.
    \item \textit{\textsc{dec-nlc}}: Decision problems with at least one nonlinear constraint.
    \item \textit{\textsc{partial-lin}}: Problems with at least one hard constraint and all constraints being linear.
    \item \textit{\textsc{soft-lin}}: Problems with all soft and linear constraints.
\end{itemize}

Table \ref{tab:intsize_distribution} shows the distribution of problem instances by the integer size. The majority of instances can be handled by \scip and \fscip. Note that \fscip was limited to handling integer sizes up to 47, whereas \scip could handle integer sizes up to 49.

\begin{table}[h]
\centering
\caption{Integer size distribution}
\begin{tabular}{@{}c@{\hspace{70pt}}r@{}} 
\toprule
integer size &\#instances \\ \midrule
0--32         & 1126 \\ 
33--47        & 77   \\
48--49         & 9    \\
50+           & 4    \\ 
\bottomrule
\end{tabular}
\label{tab:intsize_distribution}
\end{table}

\subsection{Building FiberSCIP for Pseudo-Boolean Problems}
\label{sec:fscip}

\fscip~\cite{shinano2018fiberscip} for Pseudo-Boolean Problems was built
by using {\em ug[SCIP-*,*]-libraries}~\cite{Shinano2019}.
To use the ug[SCIP-*,*]-libraries, the user has to provide glue code, which is a list of \scip plug-ins specialized for the Pseudo-Boolean problems.
Only 57 lines of glue code are necessary for the parallelization.

We used \emph{customized racing}, a parallel solving strategy where multiple \scip instances tackle the same problem simultaneously. 
Users can specify a set of parameters individually for each \scip solver thread managed by \fscip\footnote{The customized racing feature existed in 2019 already, but it was not officially released since the feature had been tested only for several customized solvers, that is \scipjack for Steiner Tree Problem, and \scipsdp for MISDP solver.}.
For the competition, we used only racing until the time limit with the following parameters:
\begin{itemize}
    \item \textit{Default}: Standard \scip settings.
    \item \textit{Aggressive Heuristics}: Prioritizes heuristics more heavily during the LP-based branch-and-bound.
    \item \textit{SAT-like Search}: Simulates a SAT-style approach.
\end{itemize}
These settings, combined with different seeds for random number generators across the 20 available cores on the competition machines, enabled each \scip instance to explore various parts of the solution space independently. 
Note that we are not processing the brand-and-bound tree of the individual runs in parallel. As far as we tested, the above strategy was more effective.

\subsection{Results}
 
\scip performed quite well, as shown in Table \ref{tab:competition_results}. Furthermore, \scip was used inside other solvers, particularly the solver that won the category \textsc{opt-lin}.

The competition was conducted with an 1 hour time limit on different machines for \scip and \fscip solvers, with the following configurations:
\begin{itemize}
    \item \scip: Intel Xeon E5-2637v4, 3.50GHz, single thread.
    \item \fscip: Intel Xeon Gold 6248, 2.50GHz, 20 threads.
\end{itemize}

 Table \ref{tab:competition_results} presents the results and rankings for\scip, \fscip, and the best of the other solvers. The displayed results are the official data as reported on the competition website  \cite{PB24}.
 
 \begin{table}[tb]
 \caption{PB24 Competition results \cite{PB24}}
\label{tab:competition_results}
\centering
\begin{tabular}{@{}l r r rrr c rr@{}}
\toprule
\multirow{2}{*}{Name} & \multirow{2}{*}{\#solvers} & \multirow{2}{*}{\#instances} & \multicolumn{3}{c}{\#best result} & & \multicolumn{2}{c}{ranking} \\
\cmidrule(lr){4-6} \cmidrule(lr){8-9}
& & & \scip & \fscip & Best Other & & \scip & \fscip \\
\midrule
\textsc{opt-lin} & 37 & 478 & 263 & 245 & \loss{279}\tablefootnote{The winner of the category \textsc{opt-lin} incorporates \scip.} & & 2 & 11 \\
\textsc{dec-lin} & 33 & 397 & 241 & 270 & \loss{312} & & 19 & 12 \\
\textsc{partial-lin} & 11 & 208 & 156 & \win{160} & 156 & & 2 & \win{1} \\
\textsc{soft-lin} & 11 & 60 & 53 & \win{55} & 47 & & 2 & \win{1} \\
\textsc{opt-nlc} & 9 & 54 & \win{37} & 36 & 30 & & \win{1} & 2 \\
\textsc{dec-nlc} & 9 & 10 & 9 & \win{10} & 9 & & 2 & \win{1} \\
\bottomrule
\end{tabular}
\end{table}
In the \textsc{opt-lin} category, \scip achieved the second place. The winner of this category, \solver{Mixed-Bag}, relied heavily on \scip. As is evident from the analysis in the next section, \scip was a major contributor to its success. In this category, the \fscip customized settings ``Aggressive Heuristics'' and ``SAT-like Search'' were less effective, as closing the dual gap proved more critical than finding the best solution. Additionally, \fscip was limited to handling integer sizes up to 47, whereas \scip could handle integer sizes up to 49. This limitation resulted in 9 instances that \fscip could not handle but \scip could, partially explaining the fewer instances solved by \fscip.

In the \textsc{dec-lin} category, \fscip outperformed \scip due to the additional settings ``Aggressive Heuristics'' and ``SAT-like Search'', which are particularly tailored for decision problems.
However, neither solver ranked high in this category, as both are inherently better suited for optimization problems.
In all other categories, both \scip and \fscip outperformed all other solvers.
Out of a total of 1,207 instances, \scip solved 759, while \fscip solved 776.

\subsection{Post-Competition Results}
\label{sec:postcompetition}

A competition is always a good opportunity to review implementations, and of course, we found further ways to improve the results afterward. 
All experiments were performed single thread on a machine with an Intel Xeon Gold 5122 CPU running at 3.60GHz. 
We reassessed the version of \scip used during the competition (\textit{comp. \scip}) and compared it with a post-competition version (\textit{post-comp. \scip}). The latter incorporates cut-based conflict analysis and employs \solver{Bliss} \cite{bliss}, rather than \solver{Nauty} \cite{nauty}, for detecting graph automorphisms to compute symmetries. 

Table \ref{tab:scip_performance_comparison} shows that post-comp. \scip solves 782 instances, compared to the 760 solved by comp. \scip. Moreover, the geometric mean of the running time over all instances (avg. time) is reduced by 9\%.

\begin{table}[h]
\centering
\caption{Post-Competition performance}
\label{tab:scip_performance_comparison}
\begin{tabular}{@{}lrrr@{}}
\toprule
Version & \#solved & avg. time [s] & quot. time  \\ \midrule
comp. \scip & 760 & 70.6 & 1.0  \\
post-comp. \scip & \win{782} & \win{63.9} & \win{0.91}  \\ \bottomrule
\end{tabular}
\end{table}

Furthermore, we analyze the individual contributions of each integrated feature by selectively disabling them. Table \ref{tab:feature_disable_impact} illustrates the impact of removing specific features on the solver's performance on instances affected by the change. 
The data underscores the critical role of symmetry handling. Disabling this feature led to a significant slowdown, with a 35\% increase in time on affected instances and 13 fewer instances solved.

Although the introduction of cut-based conflict analysis marked an improvement, its impact was relatively modest, with only 8 fewer instances solved when disabled. MIP solvers already include many forms of conflict analysis, and the cut-based conflict analysis in \scip supplements rather than fundamentally transforms the already robust conflict analysis mechanisms \cite{achterberg2007conflict,Witzig_2021}.

Given the limited number of non-linear instances (60 across \textsc{opt-nlc} and \textsc{dec-nlc} categories), the effects of disabling Flower inequalities and RLT cuts were observed on a smaller scale. However, where applicable, both features significantly speed up the solving process, as evidenced by the increase in relative time on the few affected instances.

\begin{table}[h]
\centering
\caption{Impact of disabling individual features}
\label{tab:feature_disable_impact}
\begin{tabular}{@{}lrrr@{}}
\toprule
disabled setting & \#affected instances & quot. time (on affected) & diff. \#solved \\ \midrule
Cut-based Conflict Analysis & 374 & \loss{1.02} & \loss{-8} \\
Symmetry & 213 & \loss{1.35} & \loss{-13} \\
RLT Cuts & 14 & \loss{1.26} & 0 \\ 
Flower Inequalities & 3 & \loss{1.40} & 0 \\ \bottomrule
\end{tabular}
\end{table}

\begin{figure}[ht!]
    \centering
    \caption{Running time for each solved instance}
    \includegraphics[width=1\textwidth]{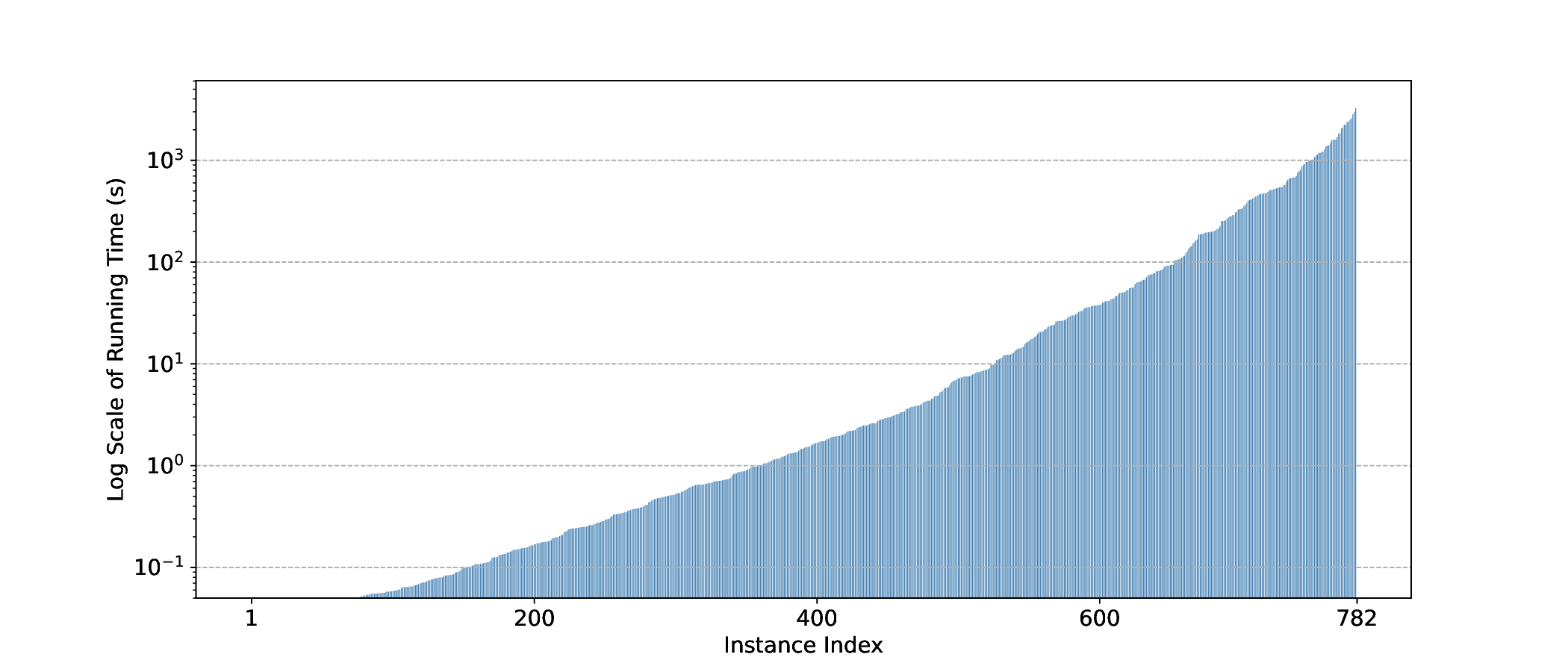}
    \label{fig:solving_time}
\end{figure}

Figure~\ref{fig:solving_time} displays the running times of \scip across all problem instances that it successfully solved. The x-axis represents the instance index, sorted in ascending order of their running times. The y-axis indicates the running times using a logarithmic scale. This distribution shows that \scip solves most instances either very quickly or not at all.
This insight was also utilized by \solver{Mixed-Bag}, the winner of the \textsc{opt-lin} category, which employed a hybrid strategy: first running  \papilo\cite{GGHpapilo} for at most 240 seconds to simplify the instance, followed by \scip on the simplified instance for up to 300 seconds. If neither  \papilo nor \scip solved the instance to optimality, a core-guided approach \cite{DevriendtGochtDemirovicNordstrom2021CuttingtotheCoreofPBO} built around \roundingsat \cite{elfers18DivideAndConquer} was employed.

% \begin{table}[h!]
% \caption{SCIP and \fscip performance with CA}
% \label{tab:competition_results_with_CA}
% \centering
% \begin{tabular}{@{}c c c c c c@{}}
% \toprule
% \multirow{2}{*}{Name} & \multirow{2}{*}{\#instances} & \multicolumn{4}{c}{\#best result} \\
% \cmidrule(lr){3-6}
% & & SCIP & SCIP+CA & \fscip & \fscip+CA \\
% \midrule
% \textsc{opt-lin}     & 478 & 263 & \loss{262} & 252 & \loss{249} \\
% \textsc{dec-lin}     & 397 & 241 & \win{255} & 273 & \win{279} \\
% \textsc{partial-lin} & 208 & 156 & \win{157} & 155 & \win{157} \\
% \textsc{soft-lin}    & 60  & 53  & 53  & 53  & \loss{51} \\
% \textsc{opt-nlc}     & 54  & 37  & 37  & 36  & 36  \\
% \textsc{dec-nlc}     & 10  & 9   & 9   & 10  & 10  \\
% \bottomrule
% \end{tabular}
% \end{table}

\section{Conclusions}
\label{sec:conclusion}
In this paper, we demonstrated the performance of \scip in the PB24 competition. Overall, solvers based on \scip won five out of the six categories. Moreover, enhancements after the competition have further improved \scip's performance as a Pseudo-Boolean solver, increasing the number of instances it can successfully solve.

% \begin{appendices}

% \section{Section title of first appendix}\label{secA1}

% \end{appendices}

\subsubsection*{Acknowledgments}
We want to thank Olivier Roussel for his efforts in organizing the competition. Part of the work for this article has been conducted in the Research Campus MODAL funded by the German Federal Ministry of Education and Research (BMBF) (fund numbers 05M14ZAM, 05M20ZBM).

\bibliographystyle{plainnat}
\bibliography{bibliography}% common bib file

\begin{thebibliography}{54}
\providecommand{\natexlab}[1]{#1}
\providecommand{\url}[1]{\texttt{#1}}
\expandafter\ifx\csname urlstyle\endcsname\relax
  \providecommand{\doi}[1]{doi: #1}\else
  \providecommand{\doi}{doi: \begingroup \urlstyle{rm}\Url}\fi

\bibitem[Aardal and Lenstra(2002)]{Aardal2002}
Karen Aardal and Arjen~K. Lenstra.
\newblock Hard equality constrained integer knapsacks.
\newblock In William~J. Cook and Andreas~S. Schulz, editors, \emph{Integer Programming and Combinatorial Optimization}, pages 350--366, Berlin, Heidelberg, 2002. Springer Berlin Heidelberg.
\newblock ISBN 978-3-540-47867-6.
\newblock \doi{10.1007/3-540-47867-1_25}.

\bibitem[Achterberg(2007{\natexlab{a}})]{Achterberg07Thesis}
Tobias Achterberg.
\newblock \emph{Constraint Integer Programming}.
\newblock Dissertation, Technische Universität Berlin, 01 2007{\natexlab{a}}.

\bibitem[Achterberg(2007{\natexlab{b}})]{achterberg2007conflict}
Tobias Achterberg.
\newblock Conflict analysis in mixed integer programming.
\newblock \emph{Discrete Optimization}, 4\penalty0 (1):\penalty0 4--20, 2007{\natexlab{b}}.
\newblock \doi{10.1016/j.disopt.2006.10.006}.

\bibitem[Adams and Sherali(1986)]{adams1986tight}
Warren~P Adams and Hanif~D Sherali.
\newblock A tight linearization and an algorithm for zero-one quadratic programming problems.
\newblock \emph{Management Science}, 32\penalty0 (10):\penalty0 1274--1290, 1986.
\newblock \doi{10.1287/mnsc.32.10.1274}.

\bibitem[Adams and Sherali(1990)]{adams1990linearization}
Warren~P Adams and Hanif~D Sherali.
\newblock Linearization strategies for a class of zero-one mixed integer programming problems.
\newblock \emph{Operations Research}, 38\penalty0 (2):\penalty0 217--226, 1990.
\newblock \doi{10.1287/opre.38.2.217}.

\bibitem[Adams and Sherali(1993)]{adams1993mixed}
Warren~P Adams and Hanif~D Sherali.
\newblock Mixed-integer bilinear programming problems.
\newblock \emph{Mathematical Programming}, 59\penalty0 (1):\penalty0 279--305, 1993.
\newblock \doi{10.1007/bf01581249}.

\bibitem[Bendotti et~al.(2021)Bendotti, Fouilhoux, and Rottner]{BendottiEtAl2021}
Pascale Bendotti, Pierre Fouilhoux, and C\'ecile Rottner.
\newblock Orbitopal fixing for the full (sub-)orbitope and application to the unit commitment problem.
\newblock \emph{Mathematical Programming}, 186:\penalty0 337--372, 2021.
\newblock \doi{10.1007/s10107-019-01457-1}.

\bibitem[Berthold et~al.(2008)Berthold, Heinz, and Pfetsch]{berthold2008solving}
Timo Berthold, Stefan Heinz, and Marc~E. Pfetsch.
\newblock Solving pseudo-{B}oolean problems with {SCIP}.
\newblock ZIB-Report 08-12, Zuse Institute Berlin, 2008.
\newblock URL \url{http://opus.kobv.de/zib/volltexte/2008/1095/}.

\bibitem[Bestuzheva et~al.(2024)Bestuzheva, Gleixner, and Achterberg]{AchterbergBestuzhevaGleixner2024}
Ksenia Bestuzheva, Ambros Gleixner, and Tobias Achterberg.
\newblock Efficient separation of {RLT} cuts for implicit and explicit bilinear terms.
\newblock \emph{Mathematical Programming B}, 2024.
\newblock \doi{10.1007/s10107-024-02104-0}.

\bibitem[Biere et~al.(2009)Biere, Heule, van Maaren, and Walsh]{Biere2009HandbookofSatisfiability}
Armin Biere, Marijn Heule, Hans van Maaren, and Toby Walsh, editors.
\newblock \emph{Handbook of Satisfiability}, volume 185 of \emph{Frontiers in Artificial Intelligence and Applications}.
\newblock IOS Press, Amsterdam, 2009.
\newblock ISBN 978-1-58603-929-5.
\newblock URL \url{http://dblp.uni-trier.de/db/series/faia/faia185.html}.

\bibitem[Bolusani et~al.(2024)Bolusani, Besan{\c{c}}on, Bestuzheva, Chmiela, Dion{\'{i}}sio, Donkiewicz, van Doornmalen, Eifler, Ghannam, Gleixner, Graczyk, Halbig, Hedtke, Hoen, Hojny, van~der Hulst, Kamp, Koch, Kofler, Lentz, Manns, Mexi, M\"{u}hmer, Pfetsch, Schl{\"o}sser, Serrano, Shinano, Turner, Vigerske, Weninger, and Xu]{SCIP9}
Suresh Bolusani, Mathieu Besan{\c{c}}on, Ksenia Bestuzheva, Antonia Chmiela, Jo{\~{a}}o Dion{\'{i}}sio, Tim Donkiewicz, Jasper van Doornmalen, Leon Eifler, Mohammed Ghannam, Ambros Gleixner, Christoph Graczyk, Katrin Halbig, Ivo Hedtke, Alexander Hoen, Christopher Hojny, Rolf van~der Hulst, Dominik Kamp, Thorsten Koch, Kevin Kofler, Jurgen Lentz, Julian Manns, Gioni Mexi, Erik M\"{u}hmer, Marc~E. Pfetsch, Franziska Schl{\"o}sser, Felipe Serrano, Yuji Shinano, Mark Turner, Stefan Vigerske, Dieter Weninger, and Lixing Xu.
\newblock {The SCIP Optimization Suite 9.0}.
\newblock Technical report, Optimization Online, February 2024.
\newblock URL \url{https://optimization-online.org/2024/02/the-scip-optimization-suite-9-0/}.

\bibitem[Chai and Kuehlmann(2005)]{ChaiKuehlmann05Galena}
D.~Chai and A.~Kuehlmann.
\newblock A fast pseudo-boolean constraint solver.
\newblock \emph{IEEE Transactions on Computer-Aided Design of Integrated Circuits and Systems}, 24\penalty0 (3):\penalty0 305--317, 2005.
\newblock \doi{10.1109/TCAD.2004.842808}.

\bibitem[Chen et~al.(2024)Chen, Lin, Hu, and Cai]{ChenLinHuCai2024ParLSPBO}
Zhihan Chen, Peng Lin, Hao Hu, and Shaowei Cai.
\newblock {ParLS-PBO: A Parallel Local Search Solver for Pseudo Boolean Optimization}.
\newblock In Paul Shaw, editor, \emph{30th International Conference on Principles and Practice of Constraint Programming (CP 2024)}, volume 307 of \emph{Leibniz International Proceedings in Informatics (LIPIcs)}, pages 5:1--5:17, Dagstuhl, Germany, 2024. Schloss Dagstuhl -- Leibniz-Zentrum f{\"u}r Informatik.
\newblock ISBN 978-3-95977-336-2.
\newblock \doi{10.4230/LIPIcs.CP.2024.5}.

\bibitem[Crama and Rodríguez-Heck(2017)]{CramaR17}
Yves Crama and Elisabeth Rodríguez-Heck.
\newblock A class of valid inequalities for multilinear 0--1 optimization problems.
\newblock \emph{Discrete Optimization}, 25:\penalty0 28--47, 2017.
\newblock ISSN 1572-5286.
\newblock \doi{10.1016/j.disopt.2017.02.001}.

\bibitem[Del~Pia and Di~Gregorio(2021)]{DelPiaD21}
Alberto Del~Pia and Silvia Di~Gregorio.
\newblock Chv{\'a}tal rank in binary polynomial optimization.
\newblock \emph{INFORMS Journal on Optimization}, 2021.
\newblock \doi{10.1287/ijoo.2019.0049}.

\bibitem[Del~Pia and Khajavirad(2017)]{DelPiaK17}
Alberto Del~Pia and Aida Khajavirad.
\newblock {A Polyhedral Study of Binary Polynomial Programs}.
\newblock \emph{Mathematics of Operations Research}, 42\penalty0 (2):\penalty0 389--410, 2017.
\newblock \doi{10.1287/moor.2016.0804}.

\bibitem[Del~Pia and Khajavirad(2018)]{DelPiaK18}
Alberto Del~Pia and Aida Khajavirad.
\newblock The multilinear polytope for acyclic hypergraphs.
\newblock \emph{SIAM Journal on Optimization}, 28\penalty0 (2):\penalty0 1049--1076, 2018.
\newblock \doi{10.1137/16M1095998}.

\bibitem[Del~Pia and Walter(2022)]{DelPiaW22}
Alberto Del~Pia and Matthias Walter.
\newblock Simple odd $\beta$-cycle inequalities for binary polynomial optimization.
\newblock In Karen Aardal and Laura Sanit{\`a}, editors, \emph{Integer Programming and Combinatorial Optimization}, pages 181--194. Springer International Publishing, 2022.
\newblock ISBN 978-3-031-06901-7.
\newblock \doi{10.1007/978-3-031-06901-7_14}.

\bibitem[Del~Pia et~al.(2020)Del~Pia, Khajavirad, and Sahinidis]{DelPiaKS20}
Alberto Del~Pia, Aida Khajavirad, and Nikolaos~V. Sahinidis.
\newblock On the impact of running intersection inequalities for globally solving polynomial optimization problems.
\newblock \emph{Mathematical Programming Computation}, 12\penalty0 (2):\penalty0 165--191, 2020.
\newblock \doi{10.1007/s12532-019-00169-z}.

\bibitem[Devriendt et~al.(2021)Devriendt, Gocht, Demirović, Nordström, and Stuckey]{DevriendtGochtDemirovicNordstrom2021CuttingtotheCoreofPBO}
Jo~Devriendt, Stephan Gocht, Emir Demirović, Jakob Nordström, and Peter Stuckey.
\newblock Cutting to the core of pseudo-{B}oolean optimization: Combining core-guided search with cutting planes reasoning.
\newblock \emph{Proceedings of the AAAI Conference on Artificial Intelligence}, 35:\penalty0 3750--3758, 05 2021.
\newblock \doi{10.1609/aaai.v35i5.16492}.

\bibitem[Dixon and Ginsberg(2002)]{DixonGinsberg02InfernceMethods}
Heidi~E. Dixon and Matthew~L. Ginsberg.
\newblock Inference methods for a pseudo-boolean satisfiability solver.
\newblock In \emph{Eighteenth National Conference on Artificial Intelligence}, pages 635--640, USA, 2002. American Association for Artificial Intelligence.
\newblock ISBN 0262511290.

\bibitem[Doornmalen and Hojny(2024)]{DoornmalenHojny2024a}
Jasper Doornmalen and Christopher Hojny.
\newblock A unified framework for symmetry handling.
\newblock \emph{Mathematical Programming}, 0\penalty0 (0):\penalty0 null, 2024.
\newblock \doi{10.1007/s10107-024-02102-2}.

\bibitem[E{\'e}n and S{\"o}rensson(2006)]{En2006TranslatingPC}
Niklas E{\'e}n and Niklas S{\"o}rensson.
\newblock Translating pseudo-boolean constraints into {SAT}.
\newblock \emph{J. Satisf. Boolean Model. Comput.}, 2:\penalty0 1--26, 2006.
\newblock \doi{10.3233/sat190014}.
\newblock URL \url{https://api.semanticscholar.org/CorpusID:4907188}.

\bibitem[Elffers and Nordström(2018)]{elfers18DivideAndConquer}
Jan Elffers and Jakob Nordström.
\newblock Divide and conquer: Towards faster pseudo-{B}oolean solving.
\newblock In \emph{Proceedings of the Twenty-Seventh International Joint Conference on Artificial Intelligence, {IJCAI-18}}, pages 1291--1299, 2018.
\newblock \doi{10.24963/ijcai.2018/180}.

\bibitem[Franco and Martin(2009)]{FrancoMartin2009HistoryofSAT}
John Franco and John Martin.
\newblock A history of satisfiability.
\newblock \emph{Frontiers in Artificial Intelligence and Applications}, 185, 01 2009.
\newblock \doi{10.3233/978-1-58603-929-5-3}.

\bibitem[Gleixner et~al.(2017)Gleixner, Eifler, Gally, Gamrath, Gemander, Gottwald, Hendel, Hojny, Koch, Miltenberger, M{\"u}ller, Pfetsch, Puchert, Rehfeldt, Schl{\"o}sser, Serrano, Shinano, Viernickel, Vigerske, Weninger, Witt, and Witzig]{SCIP5}
Ambros Gleixner, Leon Eifler, Tristan Gally, Gerald Gamrath, Patrick Gemander, Robert~Lion Gottwald, Gregor Hendel, Christopher Hojny, Thorsten Koch, Matthias Miltenberger, Benjamin M{\"u}ller, Marc~E. Pfetsch, Christian Puchert, Daniel Rehfeldt, Franziska Schl{\"o}sser, Felipe Serrano, Yuji Shinano, Jan~Merlin Viernickel, Stefan Vigerske, Dieter Weninger, Jonas~T. Witt, and Jakob Witzig.
\newblock {The SCIP Optimization Suite 5.0}.
\newblock Technical Report 17-61, ZIB, Takustr. 7, 14195 Berlin, 2017.

\bibitem[Gleixner et~al.(2023)Gleixner, Gottwald, and Hoen]{GGHpapilo}
Ambros Gleixner, Leona Gottwald, and Alexander Hoen.
\newblock {PaPILO}: A parallel presolving library for integer and linear programming with multiprecision support.
\newblock \emph{INFORMS Journal on Computing}, 2023.
\newblock \doi{10.1287/ijoc.2022.0171.cd}.
\newblock URL \url{https://github.com/INFORMSJoC/2022.0171}.

\bibitem[Granlund and the {GMP}~development team(2023)]{gmplib}
Torbjörn Granlund and the {GMP}~development team.
\newblock \emph{{GNU} {MP}: {T}he {GNU} {M}ultiple {P}recision {A}rithmetic {L}ibrary}, 6.3.0 edition, 2023.
\newblock \url{http://gmplib.org/}.

\bibitem[Hemery and Lecoutre(2006)]{HemeryLecoutre2006AbsconPseudo}
F.~Hemery and C.~Lecoutre.
\newblock {AbsCon}, 2006.
\newblock URL \url{https://www.cril.univ-artois.fr/PB06/papers/abscon2006V2.pdf}.

\bibitem[Hojny(2024)]{Hojny2024}
Christopher Hojny.
\newblock Detecting and handling reflection symmetries in mixed-integer (nonlinear) programming.
\newblock \url{https://optimization-online.org/?p=26398}, 2024.

\bibitem[Hooker(1992)]{Hooker92GeneralizedResolution}
John Hooker.
\newblock Generalized resolution for 0–1 linear inequalities.
\newblock \emph{Annals of Mathematics and Artificial Intelligence}, 6:\penalty0 271--286, 03 1992.
\newblock \doi{10.1007/BF01531033}.

\bibitem[Joshi et~al.(2015)Joshi, Martins, and Manquinho]{JoshiMartins2015OpenWbo}
Saurabh Joshi, Ruben Martins, and Vasco Manquinho.
\newblock Generalized totalizer encoding for pseudo-{B}oolean constraints.
\newblock volume 9255, 07 2015.
\newblock ISBN 978-3-319-23218-8.
\newblock \doi{10.1007/978-3-319-23219-5_15}.

\bibitem[Junttila and Kaski(2012)]{bliss}
Tommi Junttila and Petteri Kaski.
\newblock bliss: A tool for computing automorphism groups and canonical labelings of graphs.
\newblock \url{http://www.tcs.hut.fi/Software/bliss/}, 2012.

\bibitem[Kaibel and Pfetsch(2008)]{KaibelPfetsch2008}
Volker Kaibel and Marc~E. Pfetsch.
\newblock Packing and partitioning orbitopes.
\newblock \emph{Mathematical Programming}, 114\penalty0 (1):\penalty0 1--36, 2008.
\newblock ISSN 0025-5610.
\newblock \doi{10.1007/s10107-006-0081-5}.

\bibitem[Kaibel et~al.(2011)Kaibel, Peinhardt, and Pfetsch]{KaibelEtAl2011}
Volker Kaibel, Matthias Peinhardt, and Marc~E. Pfetsch.
\newblock Orbitopal fixing.
\newblock \emph{Discrete Optimization}, 8\penalty0 (4):\penalty0 595--610, 2011.
\newblock ISSN 1572-5286.
\newblock \doi{http://dx.doi.org/10.1016/j.disopt.2011.07.001}.

\bibitem[Le~Berre and Parrain(2010)]{LeBerreParrain2010Sat4J22}
Daniel Le~Berre and Anne Parrain.
\newblock The {Sat4j} library, release 2.2.
\newblock \emph{JSAT}, 7:\penalty0 59--6, 01 2010.

\bibitem[Luteberget and Sartor(2023)]{luteberget2023feasibility}
Bj{\o}rnar Luteberget and Giorgio Sartor.
\newblock {Feasibility Jump: an LP-free Lagrangian MIP heuristic}.
\newblock \emph{Mathematical Programming Computation}, 15\penalty0 (2):\penalty0 365--388, 2023.
\newblock \doi{10.1007/s12532-023-00234-8}.

\bibitem[Margot(2009)]{Margot2010}
Fran{\c{c}}ois Margot.
\newblock Symmetry in integer linear programming.
\newblock \emph{50 Years of Integer Programming 1958-2008: From the Early Years to the State-of-the-Art}, pages 647--686, 2009.
\newblock \doi{10.1007/978-3-540-68279-0_17}.

\bibitem[Marques-Silva and Sakallah(1999)]{MarquesSakallah1999Grasp}
J.P. Marques-Silva and K.A. Sakallah.
\newblock {GRASP}: a search algorithm for propositional satisfiability.
\newblock \emph{IEEE Transactions on Computers}, 48\penalty0 (5):\penalty0 506--521, 1999.
\newblock \doi{10.1109/12.769433}.

\bibitem[Martins et~al.(2014)Martins, Manquinho, and Lynce]{MartinsManquinho2014OpenWbo}
Ruben Martins, Vasco Manquinho, and In{\^e}s Lynce.
\newblock {Open-WBO}: A modular {MaxSAT} solver,.
\newblock In Carsten Sinz and Uwe Egly, editors, \emph{Theory and Applications of Satisfiability Testing -- SAT 2014}, pages 438--445, Cham, 2014. Springer International Publishing.
\newblock ISBN 978-3-319-09284-3.
\newblock \doi{10.1007/978-3-319-09284-3_33}.

\bibitem[McKay and Piperno(2014)]{nauty}
Brendan~D. McKay and Adolfo Piperno.
\newblock Practical graph isomorphism, {II}.
\newblock \emph{Journal of Symbolic Computation}, 60:\penalty0 94--112, 2014.
\newblock \doi{10.1016/j.jsc.2013.09.003}.

\bibitem[Mexi et~al.(2023)Mexi, Berthold, Gleixner, and Nordstr{\"o}m]{MexiBertholdGleixneretal.2023}
Gioni Mexi, Timo Berthold, Ambros Gleixner, and Jakob Nordstr{\"o}m.
\newblock Improving conflict analysis in {MIP} solvers by pseudo-{B}oolean reasoning.
\newblock In \emph{29th International Conference on Principles and Practice of Constraint Programming (CP 2023)}, volume 280, pages 27:1--27:19, 2023.
\newblock \doi{10.4230/LIPIcs.CP.2023.27}.

\bibitem[Mexi et~al.(2024)Mexi, Serrano, Berthold, Gleixner, and Nordström]{mexi2024cutbasedconflictanalysismixed}
Gioni Mexi, Felipe Serrano, Timo Berthold, Ambros Gleixner, and Jakob Nordström.
\newblock Cut-based conflict analysis in mixed integer programming, 2024.
\newblock URL \url{https://arxiv.org/abs/2410.15110}.

\bibitem[Ostrowski et~al.(2011)Ostrowski, Linderoth, Rossi, and Smriglio]{OstrowskiEtAl2011}
James Ostrowski, Jeff Linderoth, Fabrizio Rossi, and Stefano Smriglio.
\newblock Orbital branching.
\newblock \emph{Mathematical Programming}, 126\penalty0 (1):\penalty0 147--178, 2011.
\newblock ISSN 0025-5610.
\newblock \doi{10.1007/s10107-009-0273-x}.

\bibitem[Roussel(2009)]{PB09}
Olivier Roussel.
\newblock Pseudo-{B}oolean competition 2009, 2009.
\newblock URL \url{http://www.cril.univ-artois.fr/PB09/}.

\bibitem[Roussel(2010)]{PB10}
Olivier Roussel.
\newblock Pseudo-{B}oolean competition 2010, 2010.
\newblock URL \url{http://www.cril.univ-artois.fr/PB10/}.

\bibitem[Roussel(2011)]{PB11}
Olivier Roussel.
\newblock Pseudo-{B}oolean competition 2011, 2011.
\newblock URL \url{http://www.cril.univ-artois.fr/PB11/}.

\bibitem[Roussel(2012)]{PB12}
Olivier Roussel.
\newblock Pseudo-{B}oolean competition 2012, 2012.
\newblock URL \url{http://www.cril.univ-artois.fr/PB12/}.

\bibitem[Roussel(2024)]{PB24}
Olivier Roussel.
\newblock Pseudo-{B}oolean competition 2024, 2024.
\newblock URL \url{http://www.cril.univ-artois.fr/PB24/}.

\bibitem[Sakai and Nabeshima(2015)]{SakaiNabeshima2015Construction}
Masahiko Sakai and Hidetomo Nabeshima.
\newblock Construction of an {ROBDD} for a {PB}-constraint in band form and related techniques for {PB}-solvers.
\newblock \emph{IEICE Transactions on Information and Systems}, E98.D:\penalty0 1121--1127, 06 2015.
\newblock \doi{10.1587/transinf.2014FOP0007}.

\bibitem[Sheini and Sakallah(2006)]{SheiniSakallah06Pueblo}
Hossein Sheini and Karem Sakallah.
\newblock Pueblo: A hybrid pseudo-{B}oolean {SAT} solver.
\newblock \emph{JSAT}, 2:\penalty0 165--189, 03 2006.
\newblock \doi{10.3233/SAT190020}.

\bibitem[Shinano et~al.(2018)Shinano, Heinz, Vigerske, and Winkler]{shinano2018fiberscip}
Yuji Shinano, Stefan Heinz, Stefan Vigerske, and Michael Winkler.
\newblock Fiber{SCIP}—a shared memory parallelization of {SCIP}.
\newblock \emph{INFORMS Journal on Computing}, 30\penalty0 (1):\penalty0 11--30, 2018.
\newblock \doi{10.1287/ijoc.2017.0762}.

\bibitem[Shinano et~al.(2019)Shinano, Rehfeldt, and Gally]{Shinano2019}
Yuji Shinano, Daniel Rehfeldt, and Tristan Gally.
\newblock An easy way to build parallel state-of-the-art combinatorial optimization problem solvers: A computational study on solving steiner tree problems and mixed integer semidefinite programs by using ug[scip-*,*]-libraries.
\newblock In \emph{2019 IEEE International Parallel and Distributed Processing Symposium Workshops (IPDPSW)}, pages 530--541, 2019.
\newblock \doi{10.1109/IPDPSW.2019.00095}.

\bibitem[Witzig et~al.(2021)Witzig, Berthold, and Heinz]{Witzig_2021}
Jakob Witzig, Timo Berthold, and Stefan Heinz.
\newblock Computational aspects of infeasibility analysis in mixed integer programming.
\newblock \emph{Mathematical Programming Computation}, 13\penalty0 (4):\penalty0 753–785, March 2021.
\newblock ISSN 1867-2957.
\newblock \doi{10.1007/s12532-021-00202-0}.

\end{thebibliography}

\end{document}